\documentclass[12pt]{amsart}
\usepackage[T1]{fontenc}
\usepackage[latin1]{inputenc}
\usepackage{typearea}
\usepackage{geometry}
\usepackage{ulem}
\usepackage{amsmath}
\usepackage{amssymb}
\usepackage{latexsym}
\usepackage{enumerate}
\usepackage{amsthm}
\usepackage[all]{xy}
\usepackage{hhline}
\usepackage{epsf} 
\usepackage{cite}
\usepackage{verbatim}
\usepackage{mathtools}
\usepackage{comment}
\usepackage{dsfont}
\usepackage{mathrsfs}

\usepackage{color}

\newtheorem{theorem}{{\sc Theorem}}[section]

\theoremstyle{remark}
\newtheorem{remark}[theorem]{{\sc Remark}}

\theoremstyle{definition}

\newcommand{\R}{\mathbb{R} }

\newcommand{\B}{\mathcal{B}}

\newcommand{\W}{\mathcal{W}}
\newcommand{\K}{\mathcal{K}}

\newcommand{\Prob}{\mathbb{P}}
\newcommand{\E}{\mathbb{E}}

\providecommand{\abs}[1]{\lvert #1\rvert}
\providecommand{\babs}[1]{\bigl\lvert #1\bigr\rvert}
\providecommand{\Babs}[1]{\Bigl\lvert #1\Bigr\rvert}

\providecommand{\BBabs}[1]{\Biggl\lvert #1\Biggr\rvert}
\providecommand{\norm}[1]{\lVert #1\rVert}

\DeclareMathOperator{\Lip}{Lip}

\renewcommand{\phi}{\varphi}
\renewcommand{\epsilon}{\varepsilon}

\renewcommand{\rho}{\varrho}
\renewcommand{\P}{\Prob}
\renewcommand{\1}{\mathds{1}}

\begin{document}
\title[Solution to the Stein equation]{A note on the solution to the standard normal Stein equation}
\author{Christian D\"obler}
\thanks{\noindent Mathematisches Institut der Heinrich Heine Universit\"{a}t D\"usseldorf\\
Email: christian.doebler@hhu.de }

\begin{abstract}  
For $z\in\R$ we reconsider the solution $f_z$ to the standard normal Stein equation, corresponding to 
the Kolmogorov distance for distributions on the real line. More precisely, we prove the Lipschitz continuity of the related function $x\mapsto \rho_z(x):=xf_z(x)$ and explain how this observation may help to remove some unnecessary technical conditions from recent Berry-Esseen bounds in the literature around the Malliavin-Stein method. We further show that the Lipschitz constants of the functions $\rho_z$ are yet not uniformly bounded in $z$ and discuss consequences of this negative result. 
\end{abstract}

\maketitle

\section{Introduction}\label{intro}
\subsection{Motivation} 
Stein's method \cite{St72, St86} is a collection of tools and techniques for estimating the distance between the distribution of a given random variable $W$ and the typically better known distribution of another random variable $Z$. This distributional distance may be measured in various metrics, which are all induced by a certain class $\mathcal{H}$ of real-valued measurable test functions $h$ such that $\E|h(Z)|<\infty$.
While it is in principle possible to setup this method for quite arbitrary distributions living on general measurable spaces, in this work we exclusively deal with the case that $Z$ has the standard normal distribution on $(\R,\B(\R))$, denoted by $\mathcal{N}(0,1)$.   

Stein's method for $\mathcal{N}(0,1)$ is based on the following observation that is called \textit{Stein's lemma} in the literature: If $X$ has distribution $\mathcal{N}(0,1)$, then for all 
absolutely continuous functions $f$ on $\R$ such that $\E|f'(X)|<\infty$, one has that 
\begin{align}\label{steinid}
\E\bigl[f'(X)\bigr]&=\E\big[Xf(X)\big].
\end{align}
Conversely, if $X$ is a real random variable such that \eqref{steinid} holds for all continuously differentiable $f$ with bounded derivative, then $X\sim\mathcal{N}(0,1)$. \medskip\\

Thus, if $W$ is close in distribution to $Z\sim \mathcal{N}(0,1)$, then 
\[\E[f'(W)-Wf(W)]\]
 should be close to zero for a large class of functions $f$. It was Stein's fundamental idea to make this heuristic precise by linking it to the distributional distance under consideration by means of a differential equation. Thus, for a class $\mathcal{H}$ of real-valued and Borel-measurable test functions $h$ on $\R$ such that $\E|h(Z)|<\infty$, consider the distance 
\[d_\mathcal{H}(W,Z):=d_\mathcal{H}\big(\mathcal{L}(W),\mathcal{L}(Z)\big):=\sup_{h\in\mathcal{H}}\Babs{\E[h(W)]-\E[h(Z)]}.\]
Then, if $f_h$ solves the \textit{Stein equation}
\begin{align}\label{steineq}
f'(x)-xf(x)=h(x)-\E[h(Z)]
\end{align}
corresponding to $h\in\mathcal{H}$, then the above distance my be written as 
\begin{align}\label{bstid}
d_\mathcal{H}(W,Z)=\sup_{h\in\mathcal{H}}\Babs{\E\big[f_h'(W)-Wf_h(W)\big]}.
\end{align} 
It turns out that the right hand side of \eqref{bstid} may often be quite effectively bounded (even uniformly in $h\in\mathcal{H}$) by exploiting the probabilistic structure of $W$ as well as the analytical properties of the solutions $f_h$, $h\in\mathcal{H}$. As can be expected, the analytical properties of $f_h$ mirror those of the test function $h$. The \textit{canonical solution} $f_h$ to \eqref{steineq} is given by 
\begin{align}\label{steinsol}
f_h(x)&=\frac{1}{\phi(x)}\int_{-\infty}^x\big(h(t)-\E[h(Z)]\big)\phi(t)dt\notag\\
&=-\frac{1}{\phi(x)}\int_x^{\infty}\big(h(t)-\E[h(Z)]\big)\phi(t)dt,
\end{align}  
where 
\[\phi(x)=\frac{1}{\sqrt{2\pi}}e^{-x^2/2}\quad\text{and}\quad\Phi(x)=\P(Z\leq x)=\int_{-\infty}^x\phi(t)dt,\quad x\in\R,\]
denote the standard normal density and distribution function, respectively. Indeed, if $f_h$ is bounded, which is typically the case, then it is easy to see that any other solution of \eqref{steineq} grows super-exponentially fast as $x\to\pm\infty$. It easily follows from integration theory that $f_h$ as in \eqref{steinsol} is indeed differentiable at the continuity points of $h$. By convention, at all other points $x$ of $\R$, the derivative $f_h'(x)$ is defined in such a way that equation \eqref{steineq} holds pointwise.
\smallskip\\

In Stein's method for normal approximation two particular choices of the test function class $\mathcal{H}$ are particularly prominent: 
\begin{enumerate}[(i)]
\item If $\mathcal{H}$ is the class $\Lip(1)$ of Lipschitz continuous functions $h$ on $\R$ with Lipschitz constant $1$, then the corresponding distance $d_\mathcal{H}$ is called \textit{Wasserstein distance} and is customarily denoted by $d_\W$. 
\item If $\mathcal{H}$ is the class of all indicators $h_z:=\1_{(-\infty,z]}$, $z\in\R$, of real half lines, then the induced distance $d_\mathcal{H}$ is called \textit{Kolmogorov distance} and is customarily denoted by $d_\K$. Note that $d_\K$ is simply the supremum norm distance between the involved distribution functions. In this case, we write $f_z:=f_{h_z}$ in short. 
\end{enumerate}
We recall here that a function $g$ on $\R$ is Lipschitz continuous, if and only if
\[\sup_{x\not=y}\frac{\abs{g(y)-g(x)}}{\abs{y-x}}<\infty\]
and, in this case, the left hand side equals the minimal Lipschitz constant of $g$, which also coincides with the essential supremum norm $\norm{g'}_\infty$ of any measurable version of the $\lambda^1$-a.e. existing derivative $g'$ of $g$.

For $h=h_z$ the solution $f_z$ has the explicit form
\begin{align}\label{steinsolk}
f_z(x)&=\begin{cases}
\frac{(1-\Phi(z))\Phi(x)}{\phi(x)},& x\leq z\\
\frac{\Phi(z)(1-\Phi(x))}{\phi(x)},&  x> z.
\end{cases}
\end{align} 
In both cases (i) and (ii), the analytical properties of the solutions $f_h$ and $f_z$ have been intensively studied. We refer to \cite[Lemma 2.3 and Lemma 2.4]{CGS} for precise results and proofs. In particular, the function $f_z$ is Lipschitz continuous with (uniform) Lipschitz constant $1$ (see \cite[Lemma 2.3]{CGS}).\smallskip\\

In what follows we write $\rho_h$ and $\rho_z$ for the functions $x\mapsto x f_h(x)$ and $x\mapsto x f_z(x)$, respectively. 

For $z\in\R$, in Lemma 2.3 of \cite{CGS} it is shown that $\rho_z$ is an increasing function on $\R$ and, starting with the work \cite{SZ}, this property has recently been exploited in several adaptions of Stein's method in order to prove novel Berry-Esseen bounds that are more accurate than previous results in this direction. In particular, in the papers \cite{LPY} and \cite{ERTZ}, dealing with Poisson and Rademacher functionals, respectively, this monotonicity has been used in connection with the respective version of Malliavin calculus. As a consequence, in both \cite[Theorem 1.15]{LPY} and \cite[Theorem 3.1]{ERTZ} it must somehow be ensured that for the random variable $F$ considered there, the quantity $\rho_z(F)$ is contained in the domain of the respective Malliavin derivative operator $D$, whenever $F$ is. In these works this is guaranteed by imposing certain additional assumptions on $F$ that, however, are neither natural nor typically straightforward to verify. Since in both these versions of Malliavin calculus the domain of $D$ is closed under Lipschitz continuous mappings, for $\rho_z(F)$ to be in the domain of $D$ it would thus be sufficient to prove the Lipschitzness of the $\rho_z$, which is in fact the main purpose of this note.   

\subsection{Main result}
If $h$ is a Lipschitz continuous function on $\R$ with Lipschitz constant $L\in[0,\infty)$, then by \cite[Lemma 2.4]{CGS} $f_h'$ is also Lipschitz continuous with Lipschitz constant $2L$ and since $f_h$ solves \eqref{steineq} it follows that 
\begin{align*}
\babs{\rho_h(y)-\rho_h(x)}&=\babs{f_h'(y)-f_h'(x)+h(x)-h(y)}\leq\babs{f_h'(y)-f_h'(x)}+\babs{h(x)-h(y)}\\
&\leq 3L|y-x|.
\end{align*}
Thus, $\rho_h$ is Lipschitz continuous with Lipschitz constant $\norm{\rho_h'}_\infty\leq3L$. In particular, $\sup_{h\in\Lip(1)} \norm{\rho_h'}_\infty\leq3<\infty$. \smallskip\\   

The following main result of this note states 
that the functions $\rho_z$, $z\in\R$, are also Lipschitz continuous but that $\sup_{z\in\R} \norm{\rho_z'}_\infty=+\infty$.

\begin{theorem}\label{rhotheo}
For any $z\in\R$ the function $\rho_z:\R\rightarrow\R$, $x\mapsto\rho_z(x)=xf_z(x)$, is Lipschitz continuous on $\R$ with minimal Lipschitz constant 
\begin{align*}
\norm{\rho_z'}_\infty&=\frac{(1+z^2)\Phi(z)(1-\Phi(z))}{\phi(z)}+|z|\big(1-\Phi(|z|)\big).
\end{align*}
In particular, one has $\sup_{z\in\R}\norm{\rho_z'}_\infty=+\infty$.
\end{theorem}

The proof is deferred to Section \ref{proofs}.

\begin{remark}
\begin{enumerate}[(a)]
\item As a consequence of Theorem \ref{rhotheo}, condition (1.23) in \cite{LPY} can be dropped without replacement. 
\item Similarly, assumption (3.3) in \cite{ERTZ}, that $\rho_z(F)+\mathbf{1}_{\{F>z\}}$ is in the domain of $D$, can be dispensed with. Note that, as in \cite{LPY}, in Theorem 3.1 of \cite{ERTZ} it is not necessary to assume that the other random variable $\mathbf{1}_{\{F>z\}}$ appearing there is contained in the domain of $D$ either. Indeed, this term may be dealt with by resorting to the alternative integration-by-parts formula in \cite[Proposition 2.2]{KRT2} and by imposing the condition (2.14) therein, which, in general, is easier to verify than that $\mathbf{1}_{\{F>z\}}$ belongs to the domain of $D$.
\end{enumerate}
\end{remark}

We finally discuss consequences of the fact that $\sup_{z\in\R}\norm{\rho_z'}_\infty=+\infty$. The function $\rho_z$ arises in Stein's method often via the identity
\begin{align}\label{neg1}
&f_z(x+u)-f_z(x)-f_z'(x)u=\int_0^u\bigl(f_z'(x+t)-f_z'(x)\bigr)dt\\
&=\int_0^u\bigl(\rho_z(x+t)-\rho_z(x)\bigr)dt+ \int_0^u\bigl(\mathbf{1}_{(-\infty,z]}(x+t)-\mathbf{1}_{(-\infty,z]}(x)\bigr)dt\notag,
\end{align}
where $x,u\in\R$ and it has been used that $f_z$ solves the Stein equation \eqref{steineq} with $h=\mathbf{1}_{(-\infty,z]}$ pointwise. If $M(z)\in[0,\infty)$ is a Lipschitz constant of $\rho_z$, which exists by Theorem \ref{rhotheo}, then the first integral on the right hand side of \eqref{neg1} may be estimated by 
\begin{align*}
&\BBabs{\int_0^u\bigl(\rho_z(x+t)-\rho_z(x)\bigr)dt}\leq M(z)\frac{u^2}{2}.
\end{align*}
In the Malliavin-Stein framework, the variable $u$ here often stands for a quantity of the form $D_yF$, where $D_y$ is the Malliavin derivative at the point $y$, and, hence one would obtain a bound of the form $M(z)|D_yF|^2/2$ for this integral. If $M:=\sup_{z\in\R}M(z)$ were in fact finite, then the resulting bound $M|D_yF|^2/2$ would indeed be very useful in such a situation, since the term on the right hand side of \eqref{neg1} is typically multiplied by another term of the order of $|u|$. However, by Theorem \ref{rhotheo}, such simple upper bounds may only be performed for a fixed value of $z$, whereas, for uniform bounds, one thus cannot circumvent the application of further sophisticated tools like the Malliavin integration-by-parts formula, for instance.

\section{Proofs}\label{proofs}
In this section we give the proof of Theorem \ref{rhotheo}. To this end we first note that from \eqref{steinsolk} we have that 
\begin{align}\label{rhoz}
\rho_z(x)&=\begin{cases}
\frac{(1-\Phi(z))x\Phi(x)}{\phi(x)},& x\leq z\\
\frac{\Phi(z)x(1-\Phi(x))}{\phi(x)},&  x> z,
\end{cases}
\end{align}  
which immediately implies that $\rho_z$ is continuous on $\R$ and that the restrictions $\rho_z|(-\infty,z]$ and $\rho_z|[z,+\infty)$ are infinitely often differentiable. Moreover, since $\Phi(-x)=1-\Phi(x)$ and $\phi(-x)=\phi(x)$, $x\in\R$, it is easy to see that 
\begin{align}\label{symrho}
\rho_{-z}(-x)&=-\rho_z(x),\quad z,x\in\R.
\end{align}

\begin{proof}[Proof of Theorem \ref{rhotheo}]
Since $\rho_z$ is continuous, to prove that $\rho_z$ is Lipschitz, it suffices to show that the restrictions $\rho_z|(-\infty,z]$ and $\rho_z|[z,+\infty)$ are both Lipschitz. Thanks to \eqref{symrho} it suffices to show this for $\rho_z|(-\infty,z]$. As the restriction of $\rho_z$ to $I_z:=(-\infty,z]$ is continuously differentiable, by the mean value theorem it suffices to prove that its derivative $\rho_z'$ is bounded on $I_z$. For $x\in\R$ we let $H(x):=\Phi(x)/\phi(x)$. Then, 
\[H'(x)=\frac{\phi(x)+x\Phi(x)}{\phi(x)}=1+xH(x),\quad x\in\R,\]
and one computes 
\begin{align}\label{p1}
\rho_z'(x)&=f_z(x)+xf_z'(x)=(1-\Phi(z))\bigl(H(x)+xH'(x)\bigr)\notag\\
&=(1-\Phi(z))\big(x+(1+x^2)H(x)\big)=(1-\Phi(z))\frac{x\phi(x)+(1+x^2)\Phi(x)}{\phi(x)}\notag\\
&=:(1-\Phi(z))G(x),
\end{align}  
where 
\[G(x):=x+(1+x^2)H(x)=\frac{x\phi(x)+(1+x^2)\Phi(x)}{\phi(x)}=:\frac{N(x)}{\phi(x)},\quad x\in\R.\]
Now
\begin{align}\label{p4}
N'(x)&=\phi(x)-x^2\phi(x)+(1+x^2)\phi(x)+2x\Phi(x)=2T(x),
\end{align}
where 
\[T(x):=\phi(x)+x\Phi(x).\]
We have $\lim_{x\to-\infty}T(x)=0$ and 
\[T'(x)=-x\phi(x)+\Phi(x)+x\phi(x)=\Phi(x)>0\]
for any $x\in I_z$ so that $T$ is increasing on $I_z$ and thus positive on $I_z$. As $\lim_{x\to-\infty}N(x)=0$, by \eqref{p4} this implies that $N$ and, hence, also $G$ is positive on $I_z$, which in turn implies that $\rho_z'>0$ on $I_z$. Thus, $\rho_z$ is increasing on $I_z$ and from \eqref{symrho} it follows again that $\rho_z$ is increasing on $\R$.   
Next, 
\begin{align}\label{p2}
G'(x)&=\frac{\phi(x)\bigl(\phi(x)-x^2\phi(x)+2x\Phi(x)+(1+x^2)\phi(x) \bigr) +x\phi(x)\big(x\phi(x)+(1+x^2)\Phi(x) \big)  }{\phi(x)^2}\notag\\
&=\frac{2\phi(x)+2x\Phi(x)+x^2\phi(x)+x(1+x^2)\Phi(x)}{\phi(x)}\notag\\
&=\frac{(2+x^2)\phi(x)+(x^3+3x)\Phi(x)}{\phi(x)}=:\frac{U(x)}{\phi(x)}.
\end{align}
If we can show that $U(x)>0$ for each $x\in I_z$, then $\rho_z''>0$ on $I_z$ in such a way that $\rho_z'$ is increasing on $I_z$ and hence
\begin{align}\label{p3}
\sup_{x\in I_z}\babs{\rho_z'(x)}&=\sup_{x\in I_z}\rho_z'(x)=\rho_z'(z)=(1-\Phi(z))\frac{z\phi(z)+(1+z^2)\Phi(z)}{\phi(z)}\notag\\
&=\frac{(1+z^2)\Phi(z)(1-\Phi(z))}{\phi(z)}+z(1-\Phi(z))=:L(z)<\infty.
\end{align} 
Now, 
\begin{align}\label{p5}
U'(x)&=2x\phi(x)-x(2+x^2)\phi(x)+(3x^2+3)\Phi(x)+(x^3+3x)\phi(x)\notag\\
&=3\big(x\phi(x)+(x^2+1)\Phi(x)\big)=:3S(x).
\end{align}
We have
\begin{align}\label{p6}
S'(x)&=\phi(x)-x^2\phi(x)+2x\Phi(x)+(x^2+1)\phi(x)\notag\\
&=2\big(\phi(x)+x\Phi(x)\big)=2 T(x)>0,
\end{align}
where $T(x)$ was defined and shown to be positive above.
As $\lim_{x\to-\infty}S(x)=0$, from \eqref{p6} we thus have $S>0$ on $I_z$ and, finally, as $\lim_{x\to-\infty}U(x)=0$, from \eqref{p5} we conclude that $U>0$ on $I_z$. Thus, \eqref{p3} is proved and $\rho_z$ is indeed Lipschitz continuous on $I_z$ with minimal Lipschitz constant $L(z)$.  

Thanks to \eqref{symrho}, $\rho_z|[z,\infty)$ is also Lipschitz continuous with minimal Lipschitz constant  
\begin{align*}
L(-z)&=\frac{(1+(-z)^2)\Phi(-z)(1-\Phi(-z))}{\phi(-z)}-z(1-\Phi(-z))\\
&=\frac{(1+z^2)\Phi(z)(1-\Phi(z))}{\phi(z)}-z\Phi(z). 
\end{align*}
Thus, the minimal Lipschitz constant of $\rho_z$ on $\R$ is given by 
\begin{align*}
M(z)&:=\max\big(L(z),L(-z)\big)=\frac{(1+z^2)\Phi(z)(1-\Phi(z))}{\phi(z)}+|z|\big(1-\Phi(|z|)\big),
\end{align*}
as claimed. The final claim now follows from the Gaussian Mills ratio inequality
\[\frac{1-\Phi(z)}{\phi(z)}\geq \frac{z}{1+z^2},\quad z\geq0,\]
as, for $z\geq0$, 
\begin{align*}
M(z)&=L(z)=\frac{(1+z^2)\Phi(z)(1-\Phi(z))}{\phi(z)}+z(1-\Phi(z))\\
&\geq \frac{(1+z^2)\Phi(z)(1-\Phi(z))}{\phi(z)}\geq z\Phi(z)\stackrel{z\to\infty}{\longrightarrow}+\infty.
\end{align*}
\end{proof}

\normalem
\bibliography{steinsol}{}
\bibliographystyle{alpha}
\end{document}